\UseRawInputEncoding
\documentclass[12pt,reqno]{amsart}%
\usepackage{amsfonts}
\usepackage{amsmath}
\usepackage{stmaryrd}
\usepackage{amssymb}
\usepackage[title]{appendix}
\usepackage{graphicx}%
\providecommand{\U}[1]{\protect\rule{.1in}{.1in}}
\newtheorem{theorem}{Theorem}[section]
\theoremstyle{plain}

\newtheorem{corollary}{Corollary}[section]

\newtheorem{lemma}{Lemma}[section]

\newtheorem{remark}{Remark}

\numberwithin{equation}{section}
\begin{document}
\title[Optimal weak estimates for Riesz potentials]{Optimal weak estimates for Riesz potentials}
\author{Liang Huang}
\address{School of Science\\ Xi'an University of Posts and Telecommunications\\Xi'an, 710121, China\\ E-mail: huangliang10@163.com}
\author{Hanli Tang}
\address{Laboratory of Mathematics and Complex Systems (Ministry of Education), School of Mathematical Sciences\\Beijing Normal University\\Beijing, 100875, China\\ E-mail: hltang@bnu.edu.cn}
\thanks{The first author is supported by National Key Research and Development Program of China (Grant No. 2020YFA0712900)  and National Natural Science Foundation of China(Grant No.11701032)
\\Corresponding author: Hanli Tang at hltang@bnu.edu.cn}


\begin{abstract}
In this note we prove a sharp reverse weak estimate for Riesz potentials
$$\|I_{s}(f)\|_{L^{\frac{n}{n-s},\infty}}\geq \gamma_sv_{n}^{\frac{n-s}{n}}\|f\|_{L^1}~~\text{for}~~0<f\in {L^1(\mathbb{R}^n)},$$
where $\gamma_s=2^{-s}\pi^{-\frac{n}{2}}\frac{\Gamma(\frac{n-s}{2})}{\Gamma(\frac{s}{2})}$. We also consider the behavior of  the best constant $\mathcal{C}_{n,s}$
of weak type estimate for Riesz potentials, and we prove $\mathcal{C}_{n,s}=O(\frac{\gamma_s}{s})$ as $s\rightarrow 0$.
\end{abstract}

\maketitle

\section{introduction}
The Riesz potentials(fractional integral operators) $I_{s}$, which play an important part in Analysis, are defined by
$$I_{s}(f)(x)=\gamma_s\int_{\mathbb{R}^n}\frac{f(x-y)}{|y|^{n-s}}dy,$$
where $0<s<n$ and $\gamma_s=2^{-s}\pi^{-\frac{n}{2}}\frac{\Gamma(\frac{n-s}{2})}{\Gamma(\frac{s}{2})}$. Such operators
were first systematically investigated by M.Riesz~\cite{R}. The $(L^p,L^q)$-boundedness of
Riesz potentials were proved by G.Hardy and J.Littlewood~\cite{HL} when $n=1$ and by S.Sobolev~\cite{S} when $n>1$. The $(L^1,L^{\frac{n}{n-s},\infty})$
-boundedness were obtained by A.Zygmund~\cite{Z}. More precisely, they established the following theorem.

\vskip0.5cm
\textbf{Theorem A. }\textit{Let} $0<s<n$ \textit{and let}
 $p,q$ \textit{satisfy} $1\leq p<q<\infty$ \textit{and} $\frac{1}{p}-\frac{1}{q}=\frac{s}{n}$, \textit{then when} $p>1$,
$$\|I_{s}(f)\|_{L^q(\mathbb{R}^n)}\leq{C(n,p,s)}\|f\|_{L^{p}(\mathbb{R}^n)}.$$
 \textit{And when }$p=1$,
$$\|I_s(f)\|_{L^{\frac{n}{n-s},\infty}(\mathbb{R}^n)}=\sup_{\lambda>0}\lambda|\{x\in{\mathbb{R}^n}:|I_{s}f|>\lambda\}|^{\frac{n-s}{n}}\leq{C(n,s)\|f\|_{L^{1}(\mathbb{R}^n)}}.$$
\vskip0.5cm
The best constant in the $(L^p,L^q)$ inequality when $p=\frac{2n}{n+s}$, $q=\frac{2n}{n-s}$ was precisely calculated by E.Lieb~\cite{L}(see also \cite{FL}), and
E.Lieb and M.Loss also offered an upper bound of $C(n,p,s)$(see chapter 4 in~\cite{LL}).

Although the best constant of $(L^p,L^q)$ estimate for  Riesz potentials has been studied for decades, to the best of the authors' knowledge there
is no result about the best constant of $(L^1,L^{\frac{n}{n-s},\infty})$ estimate for  Riesz potentials. In  this paper, we will provide some estimates for
the best constant of the weak type inequality.

In \cite{T}(see multilinear case in \cite{TW}), the second author setted up the following  limiting weak-type behavior for Riesz potentials,
\begin{align*}
\lim_{\lambda\rightarrow{0}}\lambda|\{x\in{\mathbb{R}^n}:|I_{s}f|>\lambda\}|^{\frac{n-s}{n}}=\gamma_s v_{n}^{\frac{n-s}{n}}\|f\|_{L^{1}(\mathbb{R}^n)}~\text{for}~~ 0<f\in L^{1}(\mathbb{R}^n),
\end{align*}
which implies a reverse weak estimate
\begin{align}\label{sharp reverse inequality}
\|I_s(f)\|_{L^{\frac{n}{n-s},\infty}(\mathbb{R}^n)}\geq \gamma_s v_{n}^{\frac{n-s}{n}}\|f\|_{L^{1}(\mathbb{R}^n)}~\text{for}~~ 0<f\in L^{1}(\mathbb{R}^n),
\end{align}
where $v_n$ is the volume of the unit ball in $\mathbb{R}^n$. So a natural question that arises here is whether the constant $\gamma_s v_{n}^{\frac{n-s}{n}}$ is sharp? In the paper, we will give an affirmative answer.

Let $\mathcal{C}_{n,s}$ be the best constant such that the $(L^1,L^{\frac{n}{n-s},\infty})$ estimate holds for Riesz potentials,
i.e. $$\mathcal{C}_{n,s}=\sup_{f\in{L^{1}(\mathbb{R}^n)}}\frac{\|I_s(f)\|_{L^{\frac{n}{n-s},\infty}(\mathbb{R}^n)}}{\|f\|_{L^{1}(\mathbb{R}^n)}}.$$
 Then from (\ref{sharp reverse inequality}), one can directly obtain a lower bound
for $\mathcal{C}_{n,s}$,
$$\mathcal{C}_{n,s}\geq \gamma_s v_{n}^{\frac{n-s}{n}}.$$

Our another goal in this paper is to provide upper and lower bounds of $\mathcal{C}_{n,s}$ to study the behavior of $\mathcal{C}_{n,s}$ as $s\rightarrow 0$. Our approach depends on the weak
$L^{\frac{n}{n-s}}$ norm $\interleave \cdot \interleave_{L^{\frac{n}{n-s},\infty}}$ which is defined by
$$\interleave f \interleave_{L^{\frac{n}{n-s},\infty}(\mathbb{R}^n)}=\sup_{0<|E|<\infty}|E|^{-\frac{1}{r}+\frac{n-s}{n}}\left(\int_E|f|^rdx\right)^{\frac{1}{r}},~~0<r<\frac{n}{n-s}.$$
The norm $\interleave \cdot \interleave_{L^{\frac{n}{n-s},\infty}}$ is equivalent to $\|\cdot\|_{L^{\frac{n}{n-s},\infty}}$. In fact there holds(see Exercise 1.1.12 in ~\cite{G})
\begin{align}\label{equa of weak norms}
\|f\|_{L^{\frac{n}{n-s},\infty}(\mathbb{R}^n)}\leq \interleave f \interleave_{L^{\frac{n}{n-s},\infty}(\mathbb{R}^n)}\leq (\frac{n}{n-r(n-s)})^{\frac{1}{r}}\|f\|_{L^{\frac{n}{n-s},\infty}(\mathbb{R}^n)}.
\end{align}

Closely related to the Riesz potentials is the centered fractional maximal function, which is defined by
$$M_{s}f(x)=\sup_{r>0}\frac{1}{|B(x,r)|^{1-\frac{s}{n}}}\int_{B(x,r)}|f(y)|dy,~~~0<s<n.$$
$M_s$ satisfy the same $(L^p,L^q)$ and $(L^1,L^{\frac{n}{n-s},\infty})$ inequality as $I_s$ does, see \cite{A} and \cite{MW}. For any positive function $f$ it is easy to see $M_s(f)\leq 1/\gamma(s)v_n^{\frac{s-n}{n}} I_s(f)$. Although the reverse inequality dose not hold in general, B.Muckenhoupt and R.Wheeden~\cite{MW} proved the two quantities are comparable in $L^p$ norm($1< p<\infty$) when $f$ is nonnegative.

Now let us state our main results. First of all we consider the weak estimate of $I_{s}(f)$ and $M_{s}(f)$  under the norm $\interleave \cdot \interleave_{L^{\frac{n}{n-s},\infty}}$. Surprisingly identities for the weak type
estimate of Riesz potentials and fractional maximal function can be established, which implies the two quantities are comparable in $L^{\frac{n}{n-s},\infty}$ (quasi)norm when $f\in L^1(\mathbb{R}^n)$ is nonnegative.

\begin{theorem}\label{theorem1.1}
Let $0<s<n$ and $f\in{L^1(\mathbb{R}^n)}$. When $1\leq r< \frac{n}{n-s}$,
$$\interleave I_{s}(f) \interleave_{L^{\frac{n}{n-s},\infty}(\mathbb{R}^n)} \leq \gamma_s v_{n}^{\frac{n-s}{n}}\left(\frac{n}{n-(n-s)r}\right)^{\frac{1}{r}}\|f\|_{L^1(\mathbb{R}^n)},$$
and
$$\interleave M_{s}(f) \interleave_{L^{\frac{n}{n-s},\infty}(\mathbb{R}^n)} = \left(\frac{n}{n-(n-s)r}\right)^{\frac{1}{r}}\|f\|_{L^1(\mathbb{R}^n)}.$$
Moreover if $0<f\in{L^1(\mathbb{R}^n)}$, then
$$\interleave I_{s}(f) \interleave_{L^{\frac{n}{n-s},\infty}(\mathbb{R}^n)}=\gamma_s v_{n}^{\frac{n-s}{n}}\left(\frac{n}{n-(n-s)r}\right)^{\frac{1}{r}}\|f\|_{L^1(\mathbb{R}^n)} . $$
\end{theorem}

\begin{remark}
In fact, from the proof one can obtain the reverse weak estimate holds when $0< r< \frac{n}{n-s}$. More precisely when $0< r< \frac{n}{n-s}$,
$$\interleave I_{s}(f) \interleave_{L^{\frac{n}{n-s},\infty}(\mathbb{R}^n)} \geq \gamma_s v_{n}^{\frac{n-s}{n}}\left(\frac{n}{n-(n-s)r}\right)^{\frac{1}{r}}\|f\|_{L^1(\mathbb{R}^n)},~~\text{if}~~0<f\in{L^1(\mathbb{R}^n)},$$
and
$$\interleave M_{s}(f) \interleave_{L^{\frac{n}{n-s},\infty}(\mathbb{R}^n)} \geq \left(\frac{n}{n-(n-s)r}\right)^{\frac{1}{r}}\|f\|_{L^1(\mathbb{R}^n)},~~\text{if}~~f\in{L^1(\mathbb{R}^n)}.$$
\end{remark}

Then we prove the following sharp reverse weak estimates for Riesz potentials.
\begin{theorem}\label{theorem1.2}
Let  $0<f\in{L^1(\mathbb{R}^n)}$, then
$$\|I_{s}(f)\|_{L^{\frac{n}{n-s},\infty}(\mathbb{R}^n)}\geq \gamma_sv_{n}^{\frac{n-s}{n}}\|f\|_{L^1(\mathbb{R}^n)}.$$
And the equality holds when
$f=(\frac{a}{b+|x-x_0|^2})^{\frac{n+s}{2}}$, where $a,b>0$ and $x_0\in{\mathbb{R}^n}$.
\end{theorem}

 As a corollary of Theorem~\ref{theorem1.1} and Theorem~\ref{theorem1.2},
we can obtain the following sharp reverse inequality.
\begin{corollary}\label{corollary1}
Let $f\in L^{1}(\mathbb{R}^{n})$, then
\begin{eqnarray*}
\|M_{s}f\|_{L^{\frac{n}{n-s},\infty}}\geq\|f\|_{L^1}.
\end{eqnarray*}
And the equality holds when $f(x)=h(|x-x_0|)$ where $h$ is a radial decreasing function.
\end{corollary}
At last we offer an upper and a lower bound for $\mathcal{C}_{n,s}$, which implies the behavior of the best constant $\mathcal{C}_{n,s}$ for small $s$ is optimal, i.e. $\mathcal{C}_{n,s}=O(\frac{\gamma_s}{s})=O(1)$
as $s\rightarrow 0$.
\begin{theorem}\label{theorem1.3}
When $n>2$ and $0<s<\frac{n-2}{4}$,
 $$\gamma_s v_n^{\frac{n-s}{n}}\frac{n-2-4s}{2s(n-2-s)}\leq \mathcal{C}_{n,s}\leq \gamma_s v_n^{\frac{n-s}{n}}\frac{n }{s}.$$
\end{theorem}
\begin{remark}
Besides using the rearrangement inequality to obtain an upper bound $\gamma_s v_n^{\frac{n-s}{n}}\frac{n }{s}$, we can take the heat-diffusion semi-group
as a tool(see the Appendix), which was used by E.Stein and  J.Str$\ddot{o}$mberg in \cite{SS} to study the $(L^1,L^{1,\infty})$ bound for centered  maximal function, to obtain
another upper bound which is equal to $O(\gamma_s v_n^{\frac{n-s}{n}}\frac{n }{s})=O(1)$
as $(s,n)\rightarrow(0,\infty)$.
\end{remark}

\section{The identity for $I_s(f)$ and $M_s(f)$ in $\interleave \cdot \interleave_{L^{\frac{n}{n-s},\infty}}$ }

In this section, we will prove Theorem~\ref{theorem1.1}. Without loss of generality let us assume $\|f\|_{L^1(\mathbb{R}^n)}=1$. Since $I_s(f)\leq I_s(|f|)$ and $r\geq 1$, using Minkowshi inequality one have for any measurable set $E$ with $|E|<\infty$,

\begin{align}\label{2.1}
& |E|^{-\frac{1}{r}+\frac{n-s}{n}}\left[\int_E|I_sf(x)|^rdx\right]^{\frac{1}{r}}
 \leq \gamma_s |E|^{-\frac{1}{r}+\frac{n-s}{n}}\int_{\mathbb{R}^n}\left[\int_E\frac{dx}{|x-y|^{(n-s)r}}\right]^{\frac{1}{r}}|f(y)|dy.
\end{align}
Then by Hardy Littlewood rearrangement inequality, there holds
\begin{align}\label{2.2}
\int_E\frac{dx}{|x-y|^{(n-s)r}}\leq \int_{E^*}\frac{dx}{|x|^{(n-s)r}}= v_{n}^{\frac{n-s}{n}r}\frac{n}{n-(n-s)r}|E|^{1-\frac{n-s}{n}r},
\end{align}
where $E^*$ is the symmetric rearrangement of the set $E$, i.e. $E^*$ is an open ball centered at the origin whose volume is $|E|$. Therefore by
(\ref{2.1}) and (\ref{2.2}) one can obtain
$$\interleave I_{s}(f) \interleave_{L^{\frac{n}{n-s},\infty}} \leq \gamma_s v_{n}^{\frac{n-s}{n}}\left(\frac{n}{n-(n-s)r}\right)^{\frac{1}{r}}\|f\|_{L^1}.$$

Next, let us prove when $0\leq f\in {L^1(\mathbb{R}^n)}$ and $0< r< \frac{n}{n-s}$,
\begin{align}~\label{reverse estimate for Is}
\interleave I_{s}(f)\interleave_{L^{\frac{n}{n-s},\infty}}\geq \gamma_s v_{n}^{\frac{n-s}{n}}\left(\frac{n}{n-(n-s)r}\right)^{\frac{1}{r}}\|f\|_{L^1} .
\end{align}
For any $\epsilon>0$, choose $R$ large enough such that $\int_{B_R(0)}f(y)dy=1-\epsilon$. Let $E=B_{lR}(0)$. Since
$$\int_{B_R(0)}\frac{f(y)}{|x-y|^{n-s}}dy\geq \int_{B_R(0)}\frac{f(y)}{(|x|+R)^{n-s}}dy=(1-\epsilon)(|x|+R)^{s-n},$$
then
\begin{align*}
& \interleave I_{s}(f)\interleave_{L^{\frac{n}{n-s},\infty}} \geq \gamma_s |E|^{-\frac{1}{r}+\frac{n-s}{n}}\left[\int_E\left(\int_{B_R(0)}\frac{f(y)}{|x-y|^{n-s}}dy\right)^{r}dx\right]^{\frac{1}{r}}\\
& \geq \gamma_s |E|^{-\frac{1}{r}+\frac{n-s}{n}}(1-\epsilon)\left[\int_{E}\frac{dx}{(|x|+R)^{(n-s)r}}\right]^{\frac{1}{r}}\\
& =\gamma_s v_n^{\frac{n-s}{n}}n^{\frac{1}{r}}(1-\epsilon)l^{-\frac{n}{r}+n-s}\left[\int_0^l\frac{t^{n-1}}{(t+1)^{(n-s)r}}dt\right]^{\frac{1}{r}}.
\end{align*}
By the fact that this inequality holds for any $l>0$, then letting $l\rightarrow\infty$, one obtain
$$\interleave I_{s}(f)\interleave_{L^{\frac{n}{n-s},\infty}}\geq \gamma_s v_{n}^{\frac{n-s}{n}}(1-\epsilon)\left(\frac{n}{n-(n-s)r}\right)^{\frac{1}{r}},$$
which implies (\ref{reverse estimate for Is}). And we finish the proof of the identity for Riesz potential.

For fractional maximum function $M_s$, since
\begin{align*}
& M_s(f)(x)\geq \frac{1}{v^{\frac{n-s}{n}}_n(|x|+R)^{n-s}}\int_{|y-x|\leq R+|x|}|f(y)|dy\\\nonumber
& \geq \frac{1}{v^{\frac{n-s}{n}}_n(|x|+R)^{n-s}}\int_{|y|\leq R}|f(y)|dy=\frac{1-\epsilon}{v^{\frac{n-s}{n}}_n(|x|+R)^{n-s}},
\end{align*}
then one can use the same method to get
$$ \interleave M_{s}(f)\interleave_{L^{\frac{n}{n-s},\infty}} \geq \left(\frac{n}{n-(n-s)r}\right)^{\frac{1}{r}}~~\text{when}~~0<r<\frac{n}{n-s}.$$
On the other hand,
$$\interleave M_{s}(f)\interleave_{L^{\frac{n}{n-s},\infty}}\leq \interleave 1/\gamma(s)v_n^{\frac{s-n}{n}} I_s(|f|)\interleave_{L^{\frac{n}{n-s},\infty}}=\left(\frac{n}{n-(n-s)r}\right)^{\frac{1}{r}}.$$
Thus one can obtain the desired identity for $M_s$.

\section{The sharp reverse weak estimate for $I_s$ and $M_s$}
In this section, first we prove the sharp reverse weak estimate for Riesz potentials $I_s$.
By (\ref{equa of weak norms}) and Theorem~\ref{theorem1.1}, there holds
$$\|I_{s}(f)\|_{L^{\frac{n}{n-s},\infty}}\geq \gamma_sv_{n}^{\frac{n-s}{n}}\|f\|_{L^1},~~0<f\in{L^{1}(\mathbb{R}^n)}.$$

Next, we will prove that the equality can be attained by the function $g(x)=(\frac{a}{b+|x-x_0|^2})^{\frac{n+s}{2}}$, where $a,b>0$ and $x_0\in{\mathbb{R}^n}$. Since the translation and dilation of $g$ do not change the
 ratio $\|I_{s}(g)\|_{L^{\frac{n}{n-s},\infty}}/\|g\|_{L^1}$, we only need to consider $g(x)=(\frac{2}{1+|x|^2})^{\frac{n+s}{2}}$. In our calculus we will use the stereographic projection, so we will introduce
 some notations about the stereographic projection here.

The inverse stereographic projection $\mathcal{S}:  \mathbb{R}^n \to \mathbb{S}^n \setminus \{S\}$, where $S = - e_{n+1}$ denotes the southpole, is given by
\begin{equation*}
\label{eq:stereo}
(\mathcal{S}(x))_i= \frac{2 x_i}{1 + |x|^2}, \quad i = 1,...,n, \quad (\mathcal{S}(x))_{n+1} = \frac{1-|x|^2}{1+|x|^2}.
\end{equation*}
Correspondingly, the stereographic projection is given by $\mathcal{S}^{-1}: \mathbb{S}^n \setminus \{S\} \to \mathbb{R}^n$,
\[ (\mathcal{S}^{-1}(\xi))_i = \frac{\xi_i}{1 + \xi_{n+1}}, \quad i = 1,...,n. \]
And the Jacobian of the (inverse) stereographic projection are
$$\mathcal{J}_{\mathcal{S}}(x)=\left(\frac{2}{1+|x|^2}\right)^{n}~~\text{and}~~\mathcal{J}_{\mathcal{S}^{-1}}(\xi)=(1+\xi_{n+1})^{-n}.$$

By changing of variables,
\begin{align}\label{norm of g}
& \|g\|_{L^1}=\int_{\mathbb{R}^n}(\frac{2}{1+|x|^2})^{\frac{n+s}{2}}dx=\int_{\mathbb{S}^n}(\frac{2}{1+|\mathcal{S}^{-1}(\xi)|^2})^{\frac{s-n}{2}}d\xi\nonumber\\
& =\int_{\mathbb{S}^n}(1+\xi_{n+1})^{\frac{s-n}{2}}d\xi=|\mathbb{S}^{n-1}|\int_{-1}^{1}(1+t)^{\frac{s-2}{2}}(1-t)^{\frac{n-2}{2}}dt\nonumber\\
& =\pi^{n/2}2^{\frac{s+n}{2}}\frac{\Gamma(s/2)}{\Gamma(\frac{s+n}{2})}.
\end{align}
Denote
$$c_{n,s}=\pi^{n/2}2^{\frac{s+n}{2}}\frac{\Gamma(s/2)}{\Gamma(\frac{s+n}{2})}.$$
Since
 $$|\mathcal{S}^{-1}(\xi)-\mathcal{S}^{-1}(\eta)|^2=\mathcal{J}_{\mathcal{S}^{-1}}(\xi)^{\frac{1}{n}}|\xi-\eta|^2\mathcal{J}_{\mathcal{S}^{-1}}(\eta)^{\frac{1}{n}},~~\text{for any}~~\xi, \eta\in{\mathbb{S}^n},$$
and
$$\int_{\mathbb{S}^n}\frac{d\eta}{|\xi-\eta|^{n-s}}=\frac{2^s\pi^{n/2}\Gamma(s/2)}{\Gamma(\frac{n+s}{2})}=\frac{c_{n,s}}{2^{\frac{n-s}{2}}}~~\text{for any}~~\eta\in{\mathbb{S}^n}~\text{(see D.4 in \cite{G})},$$
 one can obtain
\begin{align*}
& I_s(g)(x)=\gamma(s)\int_{\mathbb{R}^n}\frac{1}{|x-y|^{n-s}}(\frac{2}{1+|y|^2})^{\frac{n+s}{2}}dy\\
& =\gamma(s)\int_{\mathbb{S}^n}\frac{1}{|\mathcal{S}^{-1}(\xi)-\mathcal{S}^{-1}(\eta)|^{n-s}}(\frac{2}{1+|\mathcal{S}^{-1}(\eta)|^2})^{\frac{s-n}{2}}d\eta\\
& =\gamma(s)\int_{\mathbb{S}^n}\frac{1}{|\xi-\eta|^{n-s}|J_{\mathcal{S}^{-1}}(\xi)|^{\frac{n-s}{2n}} |J_{\mathcal{S}^{-1}}(\eta)|^{\frac{n-s}{2n}}}(\frac{2}{1+|\mathcal{S}^{-1}(\eta)|^2})^{\frac{s-n}{2}}d\eta\\
& =\gamma(s)\int_{\mathbb{S}^n}\frac{d\eta}{|\xi-\eta|^{n-s}}(1+\xi_{n+1})^{\frac{n-s}{2}}=\gamma(s)\frac{c_{n,s}}{(1+|x|^2)^{\frac{n-s}{2}}}.
\end{align*}
Thus for any $\lambda>0$,
\begin{align}\label{estimate for I_s(g)}
|\{I_s(g)>\lambda\}|=v_{n}\left((\frac{\gamma(s)c_{n,s}}{\lambda})^{\frac{2}{n-s}}-1\right)^{\frac{n}{2}}.
\end{align}
Therefore combining (\ref{norm of g}) and (\ref{estimate for I_s(g)}) one have
$$\frac{\|I_{s}(g)\|_{L^{\frac{n}{n-s},\infty}}}{\|g\|_{L^1}}=v_{n}^{\frac{n-s}{n}}\sup_{\lambda>0}\left(\gamma(s)^{\frac{2}{n-s}}-(\frac{\lambda}{c_{n,s}})^{\frac{2}{n-s}}\right)^{\frac{n-s}{2}}=\gamma(s)v_{n}^{\frac{n-s}{n}}.$$

Next let us prove the sharp reverse weak estimate for $M_s$.
By the identity in Theorem~\ref{theorem1.1} for $M_s$ and (\ref{equa of weak norms})
one can find  for any $f\in{L^1}$,
\begin{align}\label{reverse estimate for Ms}
\|M_s(f)\|_{L^{\frac{n}{n-s},\infty}}\geq \|f\|_{L^1}.
\end{align}
On the other hand, since $M_s(f)\leq 1/\gamma(s)v_n^{\frac{s-n}{n}} I_s(f)$ and we already proved that the function $g=(\frac{a}{b+|x-x_0|^2})^{\frac{n+s}{2}}$ satisfies $\|I_{s}(g)\|_{L^{\frac{n}{n-s},\infty}}=\gamma(s)v_{n}^{\frac{n-s}{n}}\|g\|_{L^1}$, then by (\ref{reverse estimate for Ms}) the following equality holds
\begin{align}\label{equality for Ms}
\|M_s(g)\|_{L^{\frac{n}{n-s},\infty}}= \|g\|_{L^1}.
\end{align}

In fact, one can prove (\ref{equality for Ms}) holds for any $L^1$ function $f(x)=h(|x-x_0|)$, where $h$ is a radial decreasing function, by using an approach from~\cite{AL}. First assume $\|f\|_{L^1}=1$. Let $\delta_{x_0}$ denote the Dirac delta mass placed at the $x_0$. It is easy to check that
$$M(\delta_{x_0})(x)=\frac{1}{|B(x,|x|)|},$$
where $M$ is the centered Hardy-Littlewood maximum function. Hence, for every $\lambda>0$, there holds
$$\lambda|\{x:M(\delta_{x_0})(x)>\lambda\}|^{\frac{n}{n-s}}=1.$$
Since $h$ is a radial decreasing function with $\|h\|_{L^1}=1$, then by the Lemma 2.1 in \cite{AL}, one have
$$M(f)(x)\leq M(\delta_{x_0})(x) \mbox{~for~every~} x\in\mathbb{R}^{n}.$$
Then for any $r>0$ and $x\in \mathbb{R}^n$,
$$\frac{1}{|B(x,r)|^{1-\frac{s}{n}}}\int_{B(x,r)}f(y)dy\leq (\frac{1}{|B(x,r)|}\int_{B(x,r)}f(y)dy\|f\|^{\frac{s}{n-s}}_{L^1})^{\frac{n-s}{n}}\leq (M(\delta_{x_0})(x))^{\frac{n-s}{n}},$$
which implies that
\begin{align}\label{estimate for Ms}
\|M_{s}f\|_{L^{\frac{n}{n-s},\infty}}\leq 1=\|f\|_{L^1}.
\end{align}
Combining this inequality with (\ref{reverse estimate for Ms}), one can obtain the desired result for $M_s$.

What is noteworthy at the end of the section is that this result is also true for centered Hardy-Littlewood maximal function. That is because using the same method one can
prove (\ref{estimate for Ms}) when $s=0$, i.e.  (\ref{estimate for Ms}) is true for centered Hardy-Littlewood maximal function. On the other hand, using the limiting weak type
behavior for maximum function in~\cite{J}, (\ref{reverse estimate for Ms}) is also true for centered Hardy-Littlewood maximal function.

\section{The upper and lower bound of $\mathcal{C}_{n,s}$ }

In this section, we will provide an upper and a lower bound for $\mathcal{C}_{n,s}$. Using Theorem~\ref{theorem1.1} and (\ref{equa of weak norms}) , we can get an upper bound
$$\mathcal{C}_{n,s}\leq \gamma_s\frac{n }{s}v_n^{\frac{n-s}{n}}.$$

To obtain the lower bound, we will use the following
formula(see section 5.10 in~\cite{LL}). Let $0<\alpha<n$, $0<s<n$ and $\alpha+s<n$, then
\begin{eqnarray}\label{equation 4.1}
\int_{\mathbb{R}^{n}}\frac{1}{|x-y|^{n-s}}\frac{1}{|y|^{n-\alpha}}dy=C_{n,\alpha,s}\frac{1}{|x|^{n-s-\alpha}}
\end{eqnarray}
with $$C_{n,\alpha,s}=\pi^{\frac{n}{2}}\frac{\Gamma(\frac{s}{2})\Gamma(\frac{\alpha}{2})\Gamma(\frac{n-s-\alpha}{2})}{\Gamma(\frac{n-s}{2})\Gamma(\frac{n-\alpha}{2})\Gamma(\frac{s+\alpha}{2})}.$$

Now assume $n-2>4s$. Choose $f(y)=\frac{1}{|y|^{n-2}}\chi_{(|y|\leq 1)}$ and let us prove
$$\|I_s f\|_{\frac{n}{n-s},\infty}\geq \gamma_s\frac{v_n^{\frac{n-s}{n}}}{s}\frac{n-2-4s}{2(n-2-s)}\|f\|_{L^1}.$$

Since $|x|\leq \frac{1}{2}$, $|y|>1$ implies $|y-x|\geq \frac{|y|}{2}$, using (\ref{equation 4.1}) with $\alpha=2$ one have
\begin{eqnarray}\label{equation 4.2}
\frac{1}{\gamma_s}I_{s}(f)(x)&=&\int_{\mathbb{R}^{n}}\frac{1}{|x-y|^{n-s}}f(y)dy=\int_{|y|\leq 1}\frac{1}{|x-y|^{n-s}}\frac{1}{|y|^{n-2}}dy\nonumber\\
&=&\int_{\mathbb{R}^{n}}\frac{1}{|x-y|^{n-s}}\frac{1}{|y|^{n-2}}dy-\int_{|y|> 1}\frac{1}{|x-y|^{n-s}}\frac{1}{|y|^{n-2}}dy\nonumber\\
&\geq&\int_{\mathbb{R}^{n}}\frac{1}{|x-y|^{n-s}}\frac{1}{|y|^{n-2}}dy-\int_{|y|> 1}\frac{2^{n-s}}{|y|^{2n-2-s}}dy\nonumber\\
&=&\frac{c}{|x|^{n-s-2}}-d ,
\end{eqnarray}
where
\begin{eqnarray*}
c=\frac{4\pi^{n/2}}{(n-s-2)\Gamma(n/2-1)s}~\text{and}~d=\frac{2^{n-s+1}\pi^{n/2}}{(n-s-2)\Gamma(n/2)}.
\end{eqnarray*}

Choose $\lambda_0=\gamma_s(2^{n-s-2}c-d)$, since $\frac{c}{d}=\frac{n-2}{s}\frac{1}{2^{n-s}}>\frac{1}{2^{n-s-2}}$, then $\lambda_0$ is positive. Thus by (\ref{equation 4.2}), there holds
\begin{eqnarray}\label{equation 4.3}
|\{I_s f>\lambda_0\}|\geq |\{|x|\leq 1/2, \frac{c}{|x|^{n-s-\alpha}}-d>\frac{\lambda_0}{\gamma_s}\}|=v_n(\frac{1}{2})^n.
\end{eqnarray}
Using the fact
$\|f\|_{L^{1}(\mathbb{R}^{n})}=\frac{\omega_{n-1}}{2}$ and (\ref{equation 4.3}) one can obtain
\begin{eqnarray*}
& \frac{\|I_\alpha f\|_{\frac{n}{n-s},\infty}}{\|f\|_{L^1}}\geq \frac{\lambda_0 |\{I_\alpha f>\lambda_0\}|^{\frac{n-s}{n}}}{\|f\|_{L^1}}=\lambda_0v^{\frac{n-s}{n}}_n\frac{\Gamma(n/2)}{2^{n-s}\pi^{n/2}}
=\gamma_s\frac{v_n^{\frac{n-s}{n}}}{s}\frac{n-2-4s}{2(n-2-s)}.
\end{eqnarray*}
So we complete the proof of Theorem~\ref{theorem1.3} .

\section*{Appendix}
In this Appendix, we give an alternative approach to prove the $(L^1,L^{\frac{n}{n-s},\infty})$ estimate for Riesz potentials, and at the same time this approach also provide an upper bound for
$\mathcal{C}_{n,s}$, which have the same behavior with $\gamma_s v_n^{(n-s)/n}n/s$
as $(s,n)\rightarrow(0,\infty)$. First, we state a lemma (see section 3 in \cite{SS}, also see the Hopf abstract maximal ergodic
theorem in \cite{DS}) about the weak estimate of the average of the heat-diffusion semi-group $T^{t}(f)=P_t \ast f$, where
$$P_t=(4\pi t)^{-\frac{n}{2}}e^{-\frac{|x|^2}{t}}.$$
\begin{lemma}\label{lamma 4.1}
For any $f\in L^1(\mathbb{R}^n)$, there holds
$$|\{x\in \mathbb{R}^n:\sup_{s>0}\frac{1}{s}\int_0^sP_tf(x)dt>\lambda\}|\leq \frac{1}{\lambda}\|f\|_{L^1(\mathbb{R}^n)}~~~,\lambda>0.$$
\end{lemma}
Now let prove the $(L^1,L^{\frac{n}{n-s},\infty})$ estimate for Riesz potentials $I_{s}(f)$, which also can be presented by the following formula related to $T^{t}(f)$,
$$I_{s}(f)(x)=\frac{1}{\Gamma(s/2)}\int^{\infty}_{0}t^{\frac{s}{2}-1}P_{t}\ast f(x)dt.$$
We divide the integral into two parts
$$\int^{\infty}_{0}t^{\frac{s}{2}-1}P_{t}\ast f(x)dt=J_{1}(f)(x)+J_{2}(f)(x),$$
where
$$J_{1}(f)(x)=\int^{R}_{0}t^{\frac{s}{2}-1}P_{t}\ast f(x)dt,$$
$$J_{2}(f)(x)=\int^{\infty}_{R}t^{\frac{s}{2}-1}P_{t}\ast f(x)dt,$$
for some $R$ to be determined later.

Denote $\mathcal{M}^{0}f(x)=\sup\limits_{r>0}\frac{1}{r}\int^{r}_{0}P_{t}\ast f(x)dt$, then we have
\begin{eqnarray}\label{A-1}
J_{1}(f)(x)&=&\sum^{\infty}_{i=1}\int^{2^{-i+1}R}_{2^{-i}R}t^{\frac{s}{2}-1}P_{t}\ast f(x)dt\nonumber\\
&\leq&\sum^{\infty}_{i=1}\int^{2^{-i+1}R}_{2^{-i}R}(2^{-i}R)^{\frac{s}{2}-1}P_{t}\ast f(x)dt\nonumber\\
&\leq&\sum^{\infty}_{i=1}(2^{-i}R)^{\frac{s}{2}-1}2^{-i+1}R\left(\frac{1}{2^{-i+1}R}\int^{2^{-i+1}R}_{0}P_{t}\ast f(x)dt\right)\nonumber\\
&\leq&2R^{\frac{s}{2}}\frac{2^{-\frac{s}{2}}}{1-2^{-\frac{s}{2}}}\mathcal{M}^{0}f(x).
\end{eqnarray}
On the other hand, by direct computation, we obtain that
\begin{eqnarray}\label{A-2}
J_{2}(f)(x)&\leq&\int^{\infty}_{R}t^{\frac{s}{2}-1}\|P_{t}\|_{L^{\infty}}\|f\|_{L^{1}}dt\nonumber\\
&\leq&\frac{2}{n-s}(4\pi)^{-\frac{n}{2}}R^{\frac{s}{2}-\frac{n}{2}}\|f\|_{L^{1}}.
\end{eqnarray}
Combining (\ref{A-1}) and (\ref{A-2}), we obtain that
\begin{eqnarray}\label{A-3}
I_{s}(f)(x)\leq\frac{1}{\Gamma(s/2)}\left(2R^{\frac{s}{2}}\frac{2^{-\frac{s}{2}}}{1-2^{-\frac{s}{2}}}\mathcal{M}^{0}f(x)+\frac{2}{n-s}(4\pi)^{-\frac{n}{2}}R^{\frac{s}{2}-\frac{n}{2}}\|f\|_{L^{1}}\right)
\end{eqnarray}
for all $R>0$. The choice of
$$R=\Big(\frac{(4\pi)^{-\frac{n}{2}}\|f\|_{L^{1}}}{\frac{s}{2^{\frac{s}{2}-1}-1}\mathcal{M}^{0}f(x)}\Big)^{\frac{2}{n}}$$
minimizes the right side of the expression in (\ref{A-3}). Thus
\begin{eqnarray}\label{A-4}
I_{s}(f)(x)&\leq&\tau_{s}(\mathcal{M}^{0}f(x))^{\frac{n-s}{n}}\|f\|^{\frac{s}{n}}_{L^{1}},
\end{eqnarray}
where
$$\tau_{s}=2(4\pi)^{-\frac{s}{2}}(2^{\frac{s}{2}}-1)^{\frac{s-n}{n}}\frac{n}{n-s}(\frac{1}{s})^{\frac{s}{n}}\frac{1}{\Gamma(s/2)}.$$
Now using lemma~\ref{lamma 4.1} one can see that
\begin{eqnarray*}
\lambda|\{I_{s}f>\lambda\}|^{\frac{n-s}{n}}&\leq&\lambda|\{\tau_{s}(\mathcal{M}^{0}f(x))^{\frac{n-s}{n}}\|f\|^{\frac{s}{n}}_{L^{1}}>\lambda\}|^{\frac{n-s}{n}}\\
&\leq&\lambda[(\frac{\tau_{s}\|f\|^{\frac{s}{n}}_{L^{1}}}{\lambda})^{\frac{n}{n-s}}\|f\|_{L^{1}}]^{\frac{n-s}{n}}\\
&\leq&\tau_{s}\|f\|_{L^{1}}.
\end{eqnarray*}
Notice that $$2^{\frac{s}{2}}-1>\frac{\ln2}{2}s\mbox{~for~}s>0,$$
thus,
\begin{eqnarray*}
\tau_{s}
&\leq&\frac{2}{\ln2}(\frac{1}{4\pi})^{-\frac{s}{2}}\frac{1}{\Gamma(\frac{s}{2}+1)}\frac{n}{n-s}.
\end{eqnarray*}
So by this approach, one can obtain that $\mathcal{C}_{n,s}\leq \frac{2}{\ln2}(\frac{1}{4\pi})^{-\frac{s}{2}}\frac{1}{\Gamma(\frac{s}{2}+1)}\frac{n}{n-s}$ and it is easy to check that
when $(s,n)\rightarrow(0,\infty)$,
$$\frac{2}{\ln2}(\frac{1}{4\pi})^{-\frac{s}{2}}\frac{1}{\Gamma(\frac{s}{2}+1)}\frac{n}{n-s}=O(\gamma_s v_n^{\frac{n-s}{n}}\frac{n }{s}).$$

\bibliographystyle{amsalpha}

\begin{thebibliography}{10}


\vskip0.5cm

\bibitem {A} R. Adams, \textit{A note on Riesz potentials}. Duke Math. J. \textbf{42} (1975), 765--778.

\bibitem {AL} J. M. Aldaz and J. P\'{e}rez L\'{a}zaro, \textit{The best constant for the centered maximal operator on radial decreasing functions}. Math. Inequal. Appl. \textbf{14} (2011), no.~1, 173--179.

\bibitem {DS} N. Dunford and J. Schwartz, \textit{Linear Operators. Part I. General Theory}. With the Assistance of W.G. Bade and R.G. Bartle,
Reprint of the 1958 original Wiley Classics Library, a Wiley-Interscience Publication, John Wiley \& Sons, Inc., New York,

\bibitem {FL} R. Frank and E. Lieb, \textit{A new, rearrangement-free proof of the Sharp Hardy-Littlewood-Sobolev
Inequality}. Spectral Theory, Function Spaces and Inequalities. Oper. Theory: Adv. Appl. \textbf{219} (2012), 55--67.

\bibitem {G} L. Grafakos, \textit{Classical Fourier Analysis }. Third edition. Graduate Texts in Mathematics, 250. Springer, New York, 2014.

\bibitem {HL} G. Hardy and J. Littlewood, \textit{Some properties of fractional integral I }. Math. Zeit.  (1927), 565--606.

\bibitem {J} P. Janakiraman, \textit{Limiting weak-type behavior for singular integral and maximal operators}. Trans. Amer. Math. Soc. \textbf{358} (2006), no.~5, 1937--1952.

\bibitem {L}E. Lieb, \textit{Sharp constants in the Hardy--Littlewood--Sobolev and related inequalities}. Ann. of Math. (2) \textbf{118} (1983), no.~2, 349--374.

\bibitem {LL}E. Lieb and M. Loss, \textit{Analysis}. Second edition. Graduate Studies in Mathematics, 14. American Mathematical Society, Providence, RI, 2001.

\bibitem {MW}B. Muckenhoupt and R. Wheeden, \textit{Weighted norm inequalities for fractional integrals}. Trans. Amer. Math. Soc. \textbf{192} (1974), 261--274.

\bibitem {R}M. Riesz, \textit{L¡¯integrale de Riemann-Liouville et le probleme de Cauchy}. Acta Math. \textbf{81} (1949), 1--222.

\bibitem {S} S. Sobolev, \textit{On a theorem in functional analysis }(in Russian). Mat. Sob. \textbf{46} (1938), 471--497.

\bibitem {SS} E. Stein and J. Stromberg, \textit{Bagavior of maximal functions in $\mathbb{R}^{n}$ for largr $n$}. Ark. Mat. \textbf{21} (1983), 259--269.

\bibitem {T} H. Tang, \textit{Limiting weak-type behavior for fractional integral operators}(Chinese). Beijing Shifan Daxue Xuebao \textbf{52} (2016), no.~1, 5--7.

\bibitem {TW} H. Tang and G.Wang \textit{Limiting weak type behavior for multilinear fractional integrals}. Nonlinear Anal. \textbf{197} (2020), 111858, 13 pp.

\bibitem {Z} A. Zygmund, \textit{On a theorem of Marcinkiewicz concerning interpolation of operators}. Jour. de Math. Pures et Apploqu\'{e}es \textbf{35} (1956), 223--248.

\end{thebibliography}

\end{document}